\documentclass[11pt]{article}

\usepackage[truedimen,margin=25truemm]{geometry}

\usepackage{amsmath,graphicx}

\usepackage{kasai_def}
\usepackage{url}

\title{Online Low-Rank Tensor Subspace Tracking from Incomplete Data\\ by CP Decomposition using Recursive Least Squares}

\author{Hiroyuki Kasai\thanks{H. Kasai is with the Graduate School of Informatics and Engineering, The University of Electro-Communications, 1-5-1 Chofugaoka, Chofu-shi, Tokyo, 182-8585, Japan (e-mail: kasai@is.uec.ac.jp, web:www.kasailab.com)}}

\begin{document}

\maketitle
\begin{abstract}
We propose an online tensor subspace tracking algorithm based on the CP decomposition exploiting the recursive least squares (RLS), dubbed OnLine Low-rank Subspace tracking  by TEnsor CP Decomposition (OLSTEC). Numerical evaluations show that the proposed OLSTEC algorithm gives faster convergence per iteration comparing with the state-of-the-art online algorithms. 
\end{abstract}

\section{Introduction}
\label{sec:intro}

The problem of tensor subspace tacking of multidimensional data, which are naturally represented by a tensor, has been studied intensively in recent years. The usual structural assumption on a tensor is that the tensor has {\it low-rank} in every mode. The popular \emph{convex relaxation}  \cite{Liu_IEEETransPAMI_2013_s, Tomioka_Latent_2011_s, Signoretto_MachineLearning_2014_s} approach minimizes the sum of the nuclear norms of the unfolding matrices of the tensor by extending the successful results in matrix completion problem \cite{Candes_FoundCompuMath_2009_s} under theoretical performance guarantees. However, due to the limited scalability towards large-scale data of convex relaxations, the \emph{fixed-rank} non-convex approach with tensor decomposition \cite{Filipovi_MultiSysSigPro_2013_s, Kressner_BIT_2014_s} has gained big attentions recently because of superior performance in practice in despite of local minima. This also comes from the success of matrix cases \cite{Boumal_NIPS_2011,Mishra_SIAMOpt_2013,Ngo_NIPS_2012_s}. Considering that the data are sequentially acquired, or the underlying low-rank structure changes over time, online subspace tracking and estimation is essential to avoid expensive repetitive computations of batch-based algorithms. 
           
With regard to matrix-based online tracking, a representative research is the projection approximation subspace tracking (PAST) \cite{Yang_IEEESP_1995_s}. GROUSE \cite{Balzano_arXiv_2010_s} recently proposes an incremental gradient descent algorithm on the Grassmannian $\mathcal{G}(d,n)$, the space of all $d$-dimensional subspace of $\mathbb{R}^n$ \cite{Edelman98a,Absil_OptAlgMatManifold_2008}. The algorithm minimizes $\ell_2$-norm cost function. GRASTA\cite{He_CVPR_2012_s} enhances robustness against outliers by exploiting $\ell_1$-norm cost function. PETRELS \cite{Chi_IEEETransSP_2013} calculates the underlying subspace via a discounted recursive process for each row of the subspace matrix in parallel. On the other hand, as for tensor-based tracking, Nion and Sidiropoulos propose an adaptive algorithm to obtain the CP (CANDECOMP/PARAFAC)  decompositions \cite{Nion_IEEETransSP_2009}. Yu et al. also propose an accelerated online tensor learning algorithm (ALTO) based on the Tucker decomposition \cite{Yu_ICML_2015_s}. However, they do not deal with missing data presence. Mardani et al. propose an online imputation algorithm based on the CP decomposition under the presence of missing data \cite{Mardani_IEEETransSP_2015}. This considers the stochastic gradient descent (SGD) for for large-scale data. However, considering the situations where the subspace changes dramatically and the processing speed is enough faster than data acquiring speed, a faster convergence algorithm per iteration to track this change is crucial. 

This paper presents a new online tensor tracking algorithm, dubbed OLSTEC, for the partially observed high-dimensional data stream corrupted by noise. We focus on the fixed-rank tensor completion algorithm with a second-order stochastic gradient descent based on the CP decomposition exploiting the recursive least squares (RLS).
The rest of paper is organized as follows. Section 2 formulates the problem of online subspace tracking and Section 3 proposes the new algorithm. Numerical evaluations are performed in Section 4, after which we conclude in Section 5.

\section{Problem Formulation}

This paper addresses the problem of low-rank tensor completion in an online manner when the rank is a priori known or estimated. Without loss of generality, we focus on 3-order tensors of which one order increases over time. In other words, we address $\mathcal{Y} \in \mathbb{R}^{L \times W \times T}$ of which 3-rd order increases infinitely. Assuming $\mathcal{Y}_{i_1, i_2, i_3}$ are only known for some indices $(i_1, i_2, i_3) \in \Omega$, where $\Omega$ is a subset of the complete set of indices $(i_1, i_2, i_3)$, a general \emph{batch-based} fixed-rank tensor completion problem is formulated as 

\begin{equation}
\begin{array}{lll}
\label{Eq:CostFunction}
\displaystyle{\min_{\mathcal{X} \in
\mathbb{R}^{L \times W \times T}} }&   
\displaystyle{\frac{1}{2}
\| \mathcal{P}_{\Omega}(\mathcal{X}) - 
\mathcal{P}_{\Omega}(\mathcal{Y}) \|^2_F} \\
{\rm subject\ to}& {\rm rank}(\mathcal{X}) = R,
\end{array}
\end{equation}
where the operator $\mathcal{P}_{\Omega}(\mathcal{X})_{i_1,i_2,i_3} = \mathcal{X}_{i_1,i_2,i_3}$ if $(i_1, i_2, i_3) \in \Omega$ and $\mathcal{P}_{\Omega}(\mathcal{X})_{i_1,i_2,i_3}  = 0$ otherwise and (with a slight abuse of notation) $\|\cdot \|_F$ is the Frobenius norm. 
${\rm rank}(\mathcal{X})$ is the rank of $\mathcal{X}$  (see \cite{Kolda_SIAMReview_2009_s} for a detailed discussion on tensor rank).
$R \ll \{L,W,T\}$ enforces a low-rank structure. Hereafter, the $t$-th slice in the third mode of $\mathcal{Y}$, i.e. $\mathcal{Y}_{:,:,t}$ and its value at $(l,w)$, i.e., $\mathcal{Y}_{l,w,t}$, are denoted as $\mat{Y}_t$ and $[\mat{Y}_t]_{l,w}$, respectively. 

The CP decomposition that we address in this paper decomposes a tensor into a sum of component rank-one tensors \cite{Kolda_SIAMReview_2009_s}, as $\mathcal{X}  \approx \sum_{r=1}^R \vec{a}_r \circ \vec{c}_r  \circ \vec{b}_r$, where $\vec{a}_r \in \mathbb{R}^L$, $\vec{b}_r \in \mathbb{R}^T$, and $\vec{c}_r \in \mathbb{R}^W$. The symbol $\circ$ represents the vector outer product.
The {\it factor matrices} refer to the combination of the vectors from the rank-one components, 
i.e., \mat{A} = $[\vec{a}_1, \vec{a}_2, \cdots \vec{a}_R] \in \mathbb{R}^{L \times R}$ and likewise for \mat{B} and \mat{C}.
It should be noted that \mat{A}, \mat{B} and \mat{C} can be also represented by row vectors, i.e., horizontal vectors, for example,   
$\mat{A} = [(\vec{a}^1)^T, \cdots, (\vec{a}^L)^T]^T$, where $\vec{a}^l \in \mathbb{R}^{R}$. 
Thus, $\mat{Y}_t=\mat{A}{\rm diag}(\vec{b}^t)\mat{C}^T=\sum_{r=1}^R \vec{b}^t(r)\vec{a}_ r \vec{c}_r^T$. 
Then, the problem (\ref{Eq:CostFunction}) is reformulated with $\ell_2$ regularizers as \cite{Mardani_IEEETransSP_2015}

\begin{equation}  
\label{eq:Batch_Problem_Definition}
\min_{\scriptsize{\mat{A},\mat{B},\mat{C}}} 
\frac{1}{2} \| \mathcal{P}_{\Omega}(\mathcal{Y}) -  \mathcal{P}_{\Omega}(\mathcal{X}) \|_F^2 + 
\mu (\| \mat{A}\|_F^2 + \| \mat{B}\|_F^2 + \| \mat{C}\|_F^2) \nonumber
\end{equation}

\begin{equation}  
 {\rm subject\ to} \ \ \ \mat{X}_\tau =  \mat{A} {\rm diag}(\vec{b}^{\tau}) \mat{C}^T   
{\rm \ for} \ \tau = 1, ..., t.
\end{equation}  
where $\mu$ is a regularizer parameter. Consequently, considering the situation where the partially observed tensor slice ${\bf \Omega}_\tau \circledast \mat{Y}_\tau$ is acquired sequentially over time, we estimate $\{ \mat{A}, \mat{B}, \mat{C}\}$ by minimizing the exponentially weighted least squares;
        \begin{eqnarray}
        	\label{Eq:Final_Problem_Definition}
        \min_{\scriptsize{\mat{A},\mat{B},\mat{C}}}  \frac{1}{2} \sum_{\tau=1}^t \lambda^{t-\tau} 
        \biggl[ {\| {\bf \Omega}_\tau  \circledast \bigl[ \mat{Y}_\tau -  \mat{A} {\rm diag}(\vec{b}^\tau) \mat{C}^T 
        \bigr] \|_F^2} 
         +\ {\bar{\mu}(\| \mat{A}\|_F^2 + \| \mat{C}\|_F^2) + \mu[\tau] \| \vec{b}^\tau \|_2^2 } 
\biggr],
        \end{eqnarray}
where $\bar{\mu}=\mu[\tau]/\sum_{\tau=1}^t \lambda^{t-\tau}$, and $0 < \lambda \leq 1$ is the so-called forgetting parameter. $\lambda=1$ case is equivalent to the batch-based problem (\ref{eq:Batch_Problem_Definition}). The symbol $\circledast$ denotes the Hadamard Product, which is the element-wise product.

\section{Proposed Tensor Tracking: OLSTEC}
The unknown variables in (\ref{Eq:Final_Problem_Definition})  are $\mat{A}, \mat{C}$, and $\vec{b}$. Since \mat{A} and \mat{C} are non-convex set, this function is non-convex. The proposed OLSTEC algorithm, as summarized by Algorithm 1, alternates between a least-square estimation of $\vec{b}[t]$ for fixed $\mat{A}[t\!-\!1]$ and $\mat{C}[t\!-\!1]$, and a second order stochastic gradient step using the RLS on $\mat{A}[t]$ and $\mat{C}[t]$ for fixed $\vec{b}[t]$. 
It should be noted that $\mat{W}[t]$ with the square bracket indicates the calculated $\mat{W}$ after performing $t$-times updates.

\subsection{Calculation of $\vec{b}[t]$}

The estimate $\vec{b}[t]$ of $\vec{b}^t$ is obtained in a closed form by least-squares by denoting  $\vec{g}_{l,w}[t]= \vec{a}^{l}[t\!\!-\!\!1] \circledast \vec{c}^{w}[t\!\!-\!\!1] \in \mathbb{R}^{R}$ 
as 
\begin{eqnarray}
\label{eq:Problem_Definition_b_timebase}
\!\min_{\scriptsize \vec{b}^t \in \mathbb{R}^{R}} 
\!\frac{1}{2} 
\Biggl[ \sum_{l=1}^L \!\!\sum_{w=1}^W \!\!
\left(
[{\bf \Omega}_t]_{l,w}  
\!\!\left(
[\mat{Y}_t]_{l,w} \!\!-\!
 (\vec{g}_{l,w}[t])^T \vec{b}^t
  \right)\! \right)^2 
\!\!+\! \mu[t] \| \vec{b}^t \|_2^2
\Biggr] \nonumber \hspace{-0.2cm}
\end{eqnarray}
Defining $F[t]$ as the inner objective to be minimized, we obtain $\vec{b}[t]$ since $\vec{b}[t]$ satisfies $\partial F[t]/\partial \vec{b}[t] = 0$ as
\begin{eqnarray}
\label{eq:Solution_b_timebase}
\vec{b}[t] & = & 
\biggl[ 
\mu[t] \mat{I}_R +  \sum_{l=1}^L  \sum_{w=1}^W {\bf \Omega}[t]_{l,w} \vec{g}_{l,w}[t]  (\vec{g}_{l,w}[t] )^T
\biggr]^{-1} 
\biggl[
\sum_{l=1}^L \sum_{w=1}^W 
{\bf \Omega}[t]_{l,w} \mat{Y}[t]_{l,w} \vec{g}_{l,w}[t]
\biggr].
\end{eqnarray}   
\vspace*{-0.3cm}

\subsection{Calculation of $\mat{A}[t]$ and $\mat{C}[t]$ based on RLS}

The calculation of $\mat{C}[t]$ uses $\mat{A}[t\!-\!1]$, and the calculation of $\mat{A}[t]$ uses $\mat{C}[t\!-\!1]$. This paper addresses a second-order stochastic gradient based on the RLS with forgetting parameters, which has been widely used in tracking of time varying parameters in many fields. Its computation is efficient since we update the estimates recursively every time new data becomes available.

As for $\mat{A}[t]$, the problem (\ref{Eq:Final_Problem_Definition}) is reformulated as
\begin{eqnarray}
\label{eq:Problem_Definition_A_timebase}
\min_{\scriptsize \mat{A} \in \mathbb{R}^{L \times R}} \frac{1}{2} \sum_{\tau}^t \lambda^{t-\tau} 
\biggl[ \| {\bf \Omega}_\tau \circledast \bigl[ \mat{Y}_\tau - 
 \mat{A} {\rm diag}(\vec{b}[\tau]) \mat{C}[\tau\!\!-\!\!1]^T 
\bigr]  \|_F^2 \biggr] + \frac{\mu[t]}{2} \| \mat{A}\|_F^2.
\end{eqnarray}

\begin{algorithm}[t]
\caption{OLSTEC algorithm}
\label{alg:algorithm}
\begin{algorithmic}[1]
\REQUIRE{ $\{ \mat{Y}_t$ and ${\bf {\Omega}}_t \}^{\infty}_{t=1}$, $\lambda$, $\mu$}
\STATE{Initialize \{$\mat{A}[0]$, $\vec{b}[0]$, $\mat{C}[0]$\}, $\mat{Y}[0]=\mat{0}$, $(\mat{RA}_l[0])^{-1}=(\mat{RC}_w[0])^{-1}=\gamma \mat{I}_{R}, \gamma > 0$.}
\FOR{$t=1,2, \cdots$} 
	\STATE{Calculate $\vec{b}[t]$ \hfill Equation (\ref{eq:Solution_b_timebase})}
	\STATE{$\mat{X}_t = \mat{A}[t\!-\!1] {\rm diag}(\vec{b}_t) (\mat{C}[t\!-\!1])^{T}$}
	\FOR{$l=1,2, \cdots, L$} 
	\STATE{Calculate $\mat{RA}_l[t]$ \hfill Equation (\ref{eq:Update_RA})}
	\STATE{Calculate $\vec{a}^l[t]$ \hfill Equation (\ref{eq:al_final})}
	\ENDFOR
	\FOR{$w=1,2, \cdots, W$} 
	\STATE{Calculate $\mat{RC}_l[t]$ \hfill Equation (\ref{eq:Update_RC})}
	\STATE{Calculate $\vec{c}^w[t]$ \hfill Equation (\ref{eq:cw_final})}
	\ENDFOR	
\ENDFOR
\RETURN $\mat{X}_t = \mat{A}[t] {\rm diag}(\vec{b}[t]) (\mat{C}[t])^{T}$
\end{algorithmic}
\end{algorithm}

The objective function in (\ref{eq:Problem_Definition_A_timebase}) decomposes into a parallel set of smaller problems, one for each row of $\mat{A}$, as
\begin{eqnarray}
	\label{eq:problem_def_am}
	\vec{a}^l[t] & \!\! =\!\! & \defargmin_{\vec{a}^l \in \mathbb{R}^{R}} \frac{1}{2}  
	\sum_{\tau=1}^t 
	\Biggl[ \lambda^{t-\tau}
	\sum_{w=1}^W
	  [{\bf \Omega}_\tau]_{l,w}\left(
	  [\mat{Y}_\tau]_{l,w} 
	  -(\vec{a}^l)^T  
	  {\rm diag} (\vec{b}[\tau])
	\vec{c}^w[\tau\!\!-\!\!1] \right) ^2 
	\Biggr]
 + \frac{\mu[t]}{2} \| \vec{a}^l\|_2^2. \nonumber 
\end{eqnarray}

Here, denoting ${\rm diag}(\vec{b}[\tau]) \vec{c}^w[\tau\!\!-\!\!1]$ as $\vec{\alpha}_w[\tau] \in \mathbb{R}^R$, $\vec{a}^l[t]$ is obtained by setting the derivative to zero as
\begin{eqnarray}	 
\label{eq:RA_a_s}
\mat{RA}_l[t]  \vec{a}^l[t]  & = &  \vec{s}_l[t],
\end{eqnarray}	
where $\mat{RA}_l[t] \in \mathbb{R}^{R\times R}$ and $\vec{s}_l[t] \in \mathbb{R}^{R}$ are defined as 
\begin{eqnarray}	
 \mat{RA}_l[t] & \!\! = \!\! &   \sum_{\tau=1}^t 
 	\biggl[
 	\sum_{w=1}^W  \lambda^{t-\tau} [{\bf \Omega}_\tau]_{l,w} 
	\vec{\alpha}_w[\tau] \vec{\alpha}_w[\tau]^T
	\biggr] + \mu[t] \mat{I}_R	\nonumber \\
	\label{eq:s_l}
	\vec{s}_l[t] & \!\! = \!\!  &  \sum_{\tau=1}^t 
	\biggl[
	\sum_{w=1}^W  \lambda^{t-\tau}  [{\bf \Omega}_\tau]_{l,w}  [\mat{Y}_\tau]_{l,w} 
	  \vec{\alpha}_w[\tau]
	 \biggr].\nonumber
\end{eqnarray}	
Here, $\mat{RA}_l[t]$ is transformed by separating $t$-th term as 
\begin{eqnarray}
\label{eq:Update_RA}
\mat{RA}_l[t]  
 &=&  \lambda \mat{RA}_l[t\!-\!1]+  \sum_{w=1}^W [{\bf \Omega}_t]_{l,w}  \vec{\alpha}_w[t]  (\vec{\alpha}_w[t])^T  + (\mu[t]  - \lambda \mu[t\!\!-\!\!1] ) \mat{I}_R.
\end{eqnarray}
Likewise, $\vec{s}_l[t]$ is obtained as 
$\vec{s}_l[t] = \lambda \vec{s}_l[t\!-\!1] + \sum_{w=1}^W [{\bf \Omega}_t]_{l,w}$\\ $ [\mat{Y}_t]_{l,w}  \vec{\alpha}_w[t]$.
Thus, from (\ref{eq:RA_a_s}), we reformulate $\mat{RA}_l[t]$ as
\begin{eqnarray}
\label{eq:RA_a}
\mat{RA}_l[t] \vec{a}^l[t]   
	 &=&  \mat{RA}_l[t] \vec{a}^l[t\!-\!1]  - (\mu[t] - \lambda \mu[t\!-\!1]) \vec{a}^l[t\!\!-\!\!1] \nonumber \\
	&&+  \sum_{w=1}^W [{\bf \Omega}_t]_{l,w} 
	\left(
	[\mat{Y}_t]_{l,w} \!-\! \vec{\alpha}_w[t]^T \vec{a}^l[t\!-\!1] \right)  
	  \vec{\alpha}_w[t].\hspace{0.5cm}\nonumber
\end{eqnarray}
Finally, $\vec{a}^l[t]$ is obtained  as
\begin{eqnarray}
	\label{eq:al_final}
	\vec{a}^l[t] 
	&=& \vec{a}^l[t\!\!-\!\!1]  - (\mu[t] - \lambda \mu[t\!\!-\!\!1]) (\mat{RA}_l[t])^{-1}\vec{a}^l[t\!\!-\!\!1]  \nonumber  \\
	 & &+ \sum_{w=1}^W [{\bf \Omega}_{t}]_{l,w} 
	 \! \left([\mat{Y}_{t}]_{l,w} \! -\!  (\vec{\alpha}_w[t])^T \vec{a}^l[t\!\!-\!\!1] \right) \! 
	(\mat{RA}_l[t])^{-1} \vec{\alpha}_w[t].
\end{eqnarray}

Similarly, $\vec{c}^w[t]$ for \mat{C}[t] can be obtained as 
\begin{eqnarray}
	\label{eq:cw_final}
	\vec{c}^w[t] & =&  \vec{c}^w[t\!\!-\!\!1]  - (\mu[t] - \lambda \mu[t\!\!-\!\!1]) (\mat{RC}_w[t])^{-1}\vec{c}^w[t\!\!-\!\!1] \ \nonumber  \\
	&&  + \sum_{l=1}^L [{\bf \Omega}_{t}]_{l,w} \! \left([\mat{Y}_{t}]_{l,w} \! -\!  \vec{\beta}_w[t] \vec{c}^w[t\!\!-\!\!1]  \right) \! 
	(\mat{RC}_w[t])^{-1} (\vec{\beta}_w[t])^T \!\!, 
\end{eqnarray}
where $\vec{\beta}_w[\tau] \in \mathbb{R}^{1 \times R}$ is ($\vec{a}^l[\tau])^T {\rm diag}(\vec{b}[\tau])$, and $\mat{RC}_w[t]$ is defined as
\begin{eqnarray}
	\label{eq:Update_RC}
	\mat{RC}_w[t] 
	& = &  \lambda \mat{RC}_w[t\!\!-\!\!1]+  \sum_{l=1}^L [{\bf \Omega}_t]_{l,w} \vec{\beta}_w[t]^T \vec{\beta}_w[t] + (\mu[t]  - \lambda \mu[t\!\!-\!\!1] ) \mat{I}_R .
\end{eqnarray}

\subsection{Complexity and memory consumption}

With respect to computational complexity per iteration, OLSTEC requires $\mathcal{O}(| {\bf {\Omega}}_t | R^2+L R^3)$ because of $\mathcal{O}(| {\bf {\Omega}}_t | R^2)$ for $\vec{b}[t]$ in (\ref{eq:Solution_b_timebase}) and $\mathcal{O}(L R^3)$ for the inversion of $\mat{RA}_l$ and $\mat{RC}_w$ in (\ref{eq:al_final}) and (\ref{eq:cw_final}), respectively. 
As for memory consumption, $\mathcal{O}((L+W)R^2)$ is required for $\mat{RA}[t]$ and $\mat{RC}[t]$, respectively.

\section{Numerical Evaluations}

\begin{figure*}[htbp]
\begin{minipage}[b]{0.50\linewidth}
\centering
\vspace{-1.0cm}
\centerline{\includegraphics[width=1.1\linewidth]{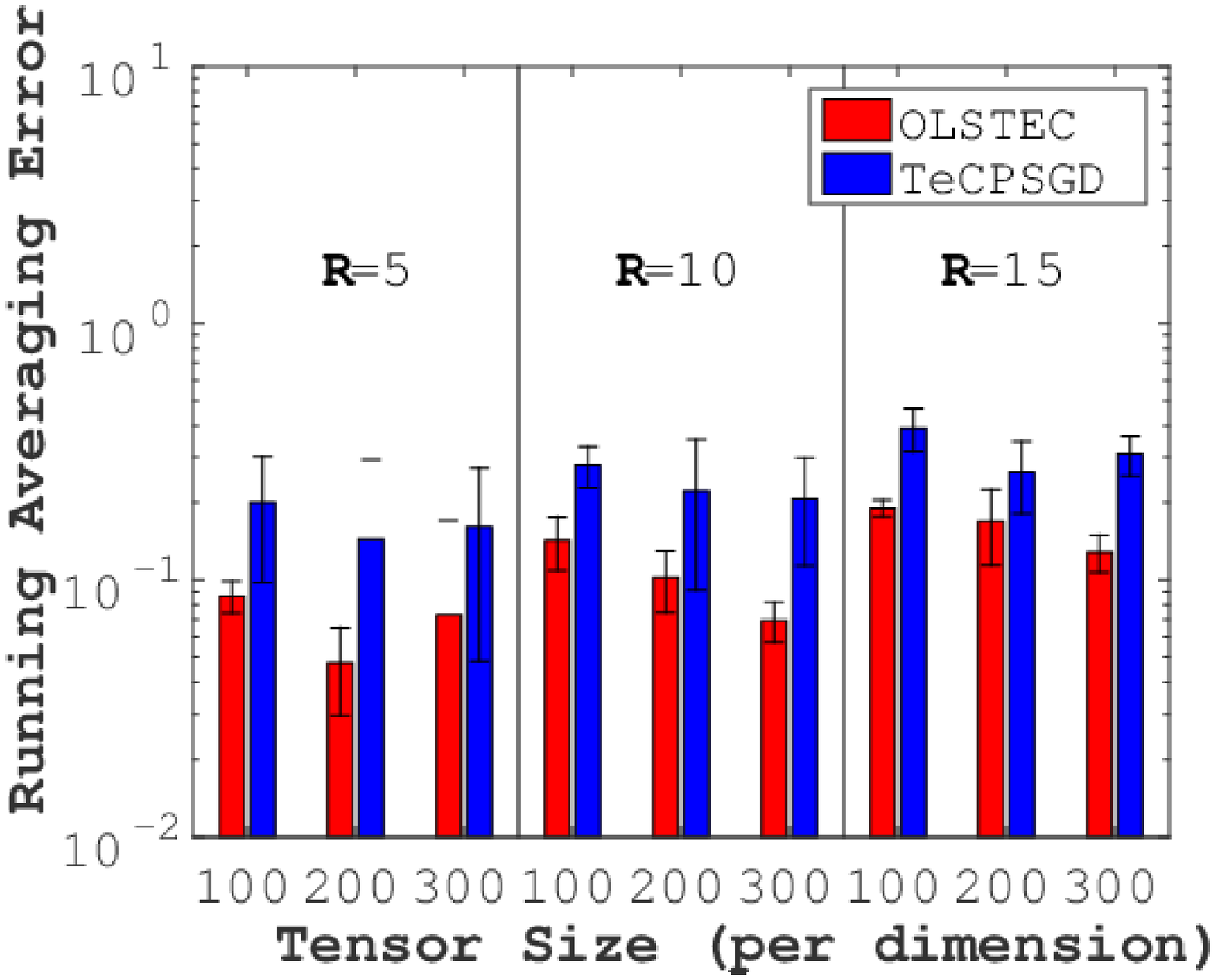}}
\vspace*{-2.5cm}
\centerline{\scriptsize (b) $\rho=0.05$}\medskip
\end{minipage}
\hspace{-0.2cm}
\begin{minipage}[b]{0.50\linewidth}
\centering
\vspace{-1.0cm}
\centerline{\includegraphics[width=1.1\linewidth]{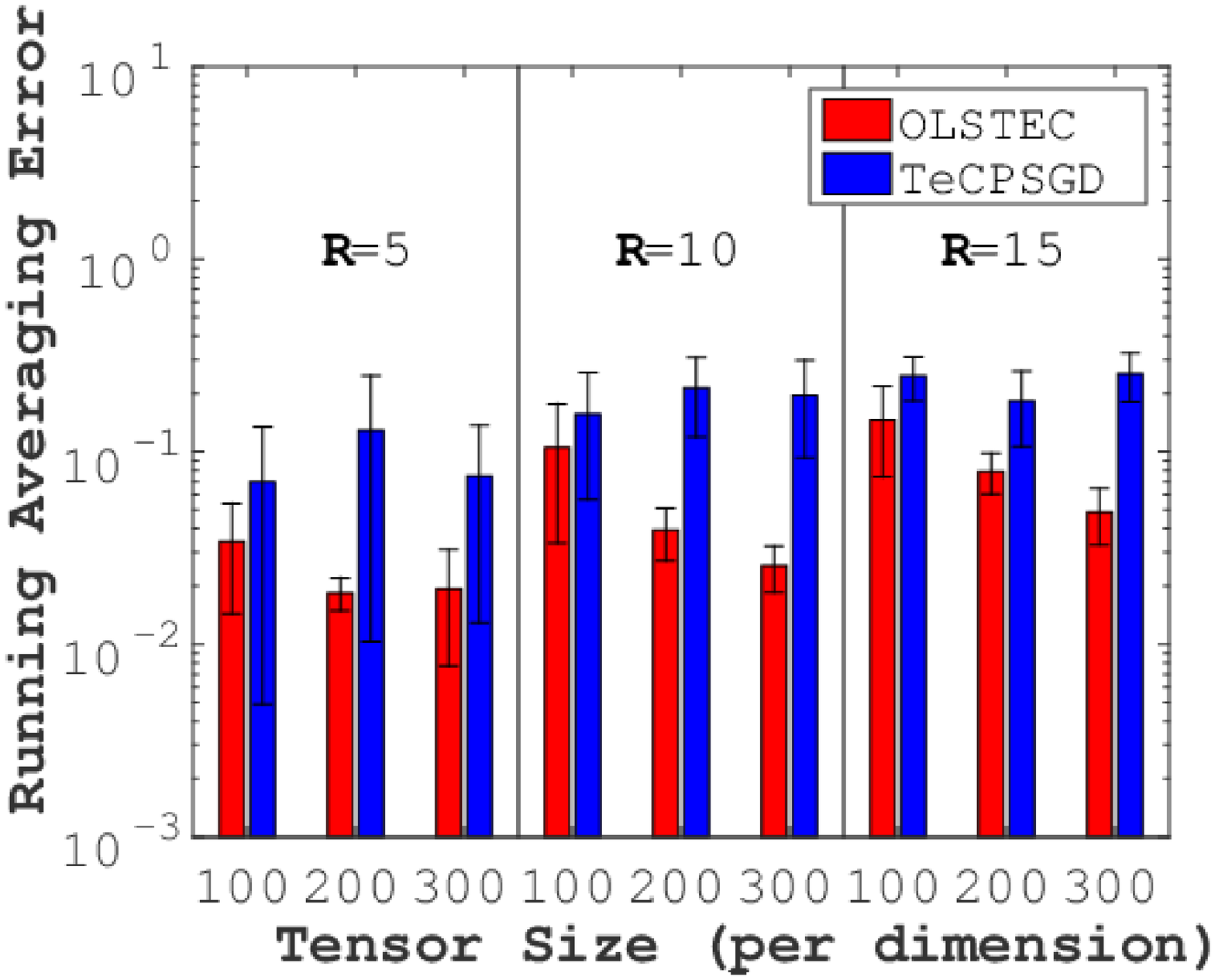}}
\vspace*{-2.5cm}
\centerline{\scriptsize (c) $\rho=0.1$}\medskip
\end{minipage}
\vspace{-0.8cm}
\caption{Running-averaging error in synthetic dataset.}
\label{fig:synthetic_data_ruuningaveraging}
\end{figure*}

\begin{figure*}[htbp]
\begin{minipage}[b]{.33\linewidth}
\centering
\centerline{\includegraphics[width=1.1\linewidth]{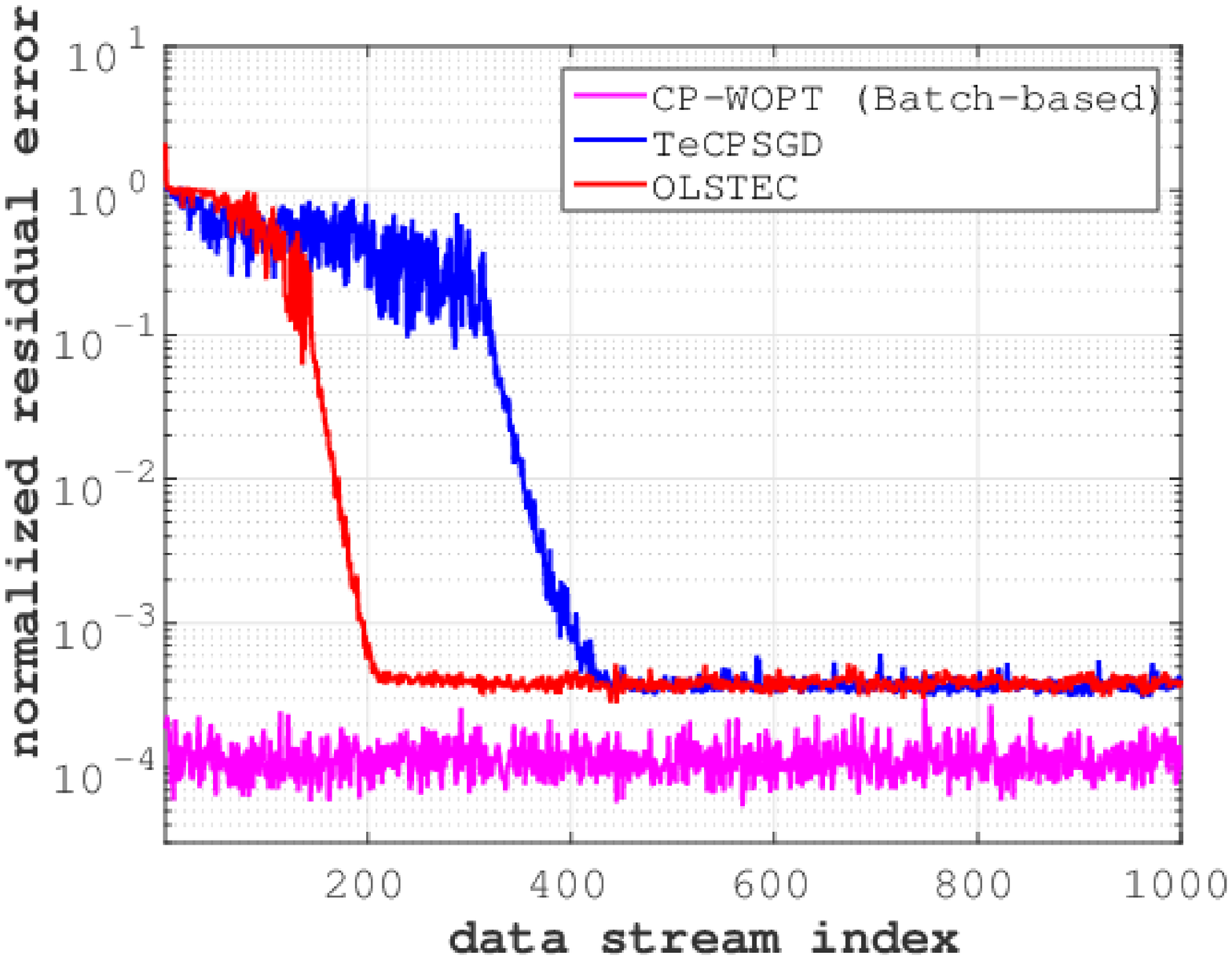}}
 \vspace*{-1.5cm}
\centerline{\scriptsize (a) Stationary subspace ($\rho = 0.1$)}\medskip
\end{minipage}
\begin{minipage}[b]{.33\linewidth}
\centering
\centerline{\includegraphics[width=1.1\linewidth]{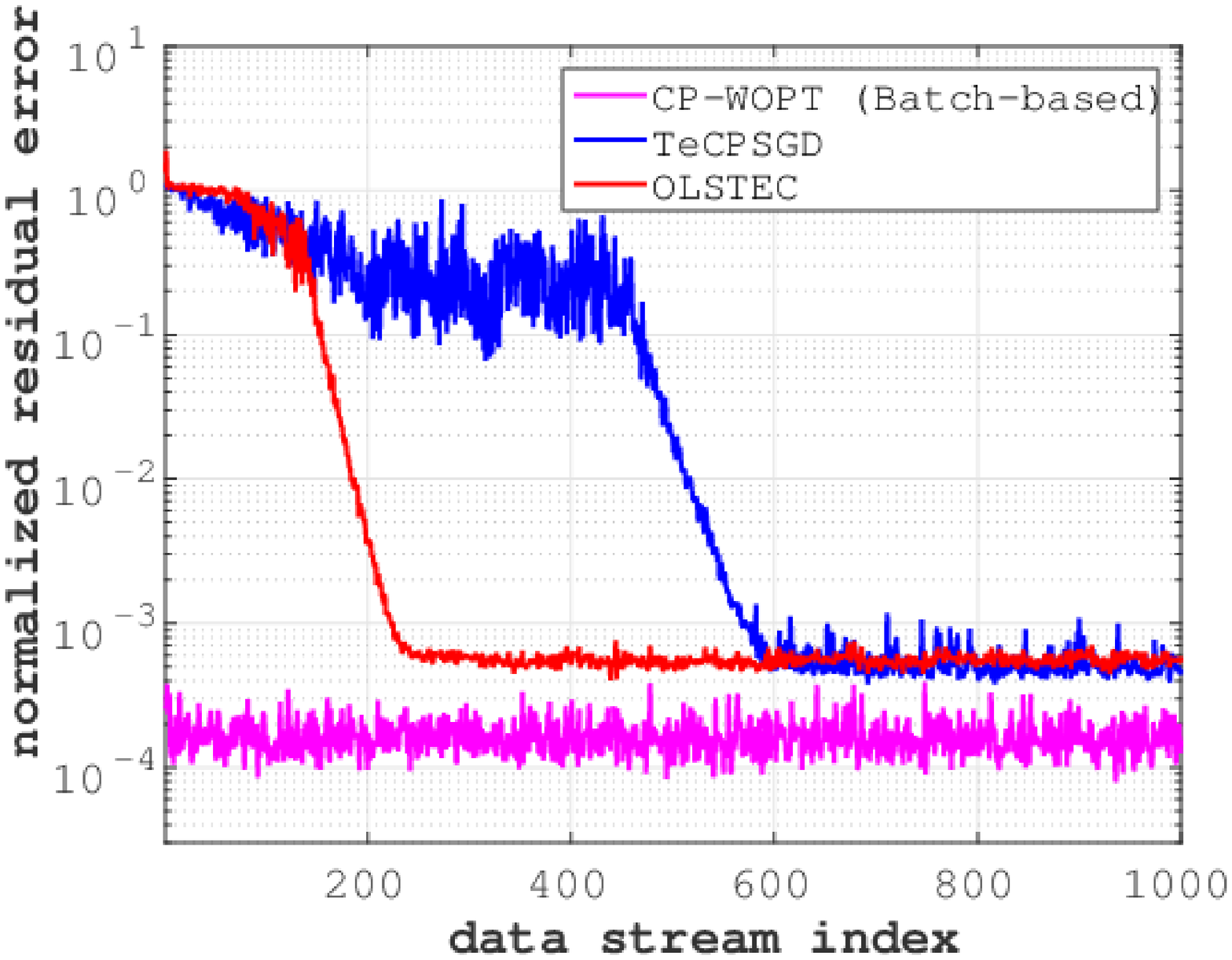}}
\vspace*{-1.5cm}
\centerline{\scriptsize (b) Stationary subspace ($\rho = 0.05$)}\medskip
\end{minipage}
\hfill
\begin{minipage}[b]{0.33\linewidth}
\centering
\centerline{\includegraphics[width=1.1\linewidth]{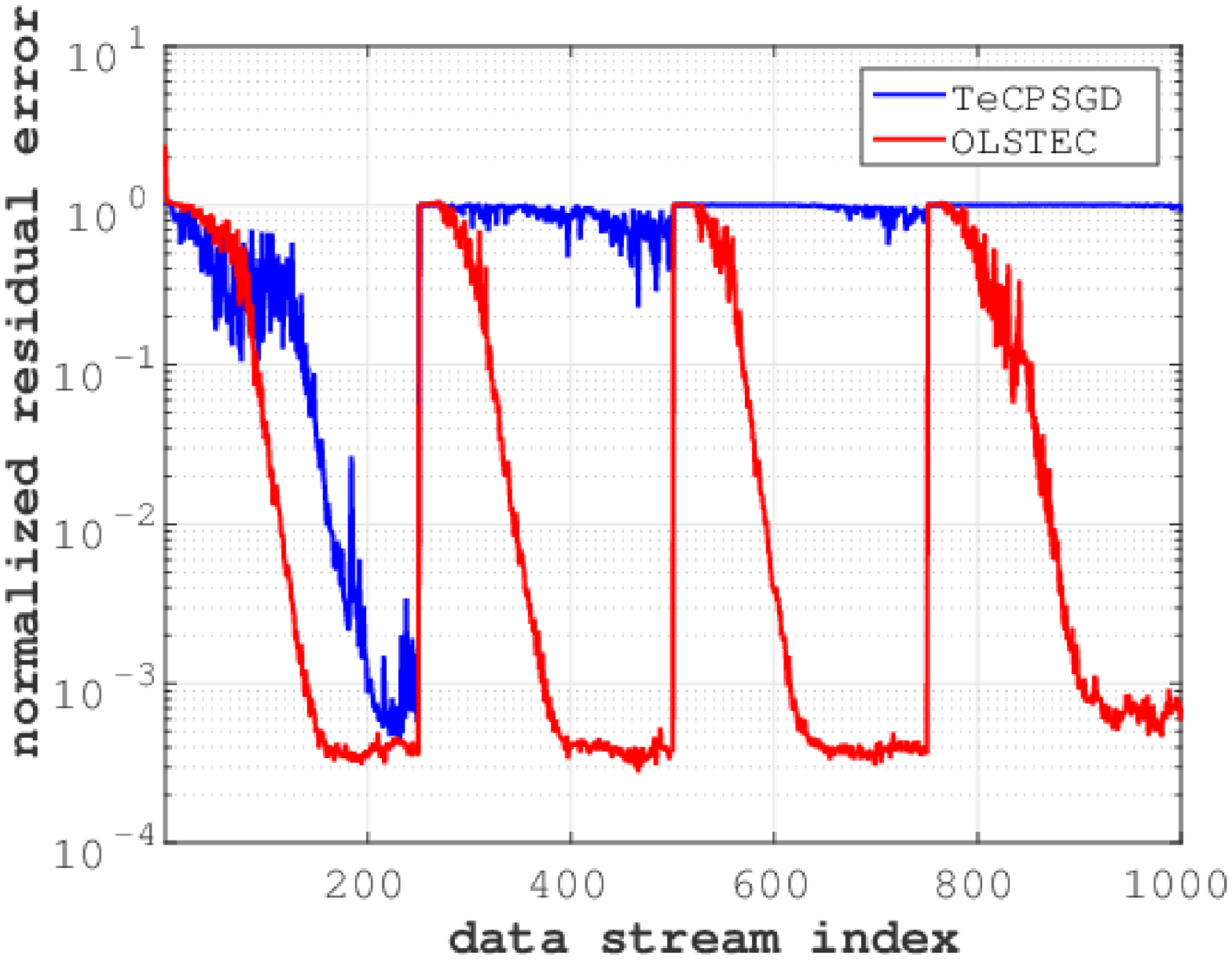}}
\vspace*{-1.5cm}
\centerline{\scriptsize (c) Dynamic subspace ($\rho = 0.1$)}\medskip
\end{minipage}
\vspace*{-0.5cm}
\caption{The normalized estimation error in synthetic dataset.}
\label{fig:synthetic_data}
\end{figure*}

We show numerical comparisons of the OLSTEC algorithm\footnote{Matlab source code is available at \url{http://www.kasailab.com/research/olstec}.} with state-of-the-art algorithms for  synthetic and real-world datasets. All the following experiments are done on a PC with $3.0$ GHz Intel Core i7 CPU and $16$ GB RAM. We first evaluates the performance of our proposed algorithm using synthetic dataset with the state-of-the-art online algorithm proposed in  \cite{Mardani_IEEETransSP_2015}, termed as ``TeCPSGD" algorithm in this paper. We first generate a low $R$-rank tensor $\mathcal{Y} \in \mathbb{R}^{L\times W \times T}$ where its factor matrices are generated with i.i.d standard Gaussian $\mathcal{N}(0, 1)$ entries, and Gaussian noise with i.i.d $\mathcal{N}(0, \epsilon^2)$ entries are added. We set $L=W=\{100,200,300\}$, $T=1000$, $R=\{5,10,15\}$, and the noise level $\epsilon=10^{-3}$. The observation ratio, $\rho$, is $\{0.1, 0.05\}$. $\mu[t]=10^{-9}$ and $\lambda=0.88$ are configured in the proposed algorithm. It should be noted that we implement TeCPSGD with our configured parameters  because the source code of TeCPSGD is not available. 
Figure \ref{fig:synthetic_data_ruuningaveraging} shows {\it the running-averaging estimation error} $\frac{1}{T} \sum_{\tau=1}^T \| \mat{X}_\tau - \mat{Y}_\tau\|^2_F/ \| \mat{Y}_\tau \|^2_F$ for each observation ratio $\rho$, where five runs are performed independently, and the results show the average with standard deviations. From these results, the proposed OLSTEC algorithm shows much lower estimation error, especially when observation ratios are lower. In addition, the standard derivations are also smaller, thus, the convergence property of the proposed algorithm is stabler than that of TeCPSGD.  Figure \ref{fig:synthetic_data} (a) and (b) show {\it the normalized residual error} $\| \mat{X}_t - \mat{Y}_t\|^2_F/ \| \mat{Y}_t \|^2_F$ when the observation ratios are $0.1$ and $0.05$, respectively. Additionally, we show, as reference, the result of CP-WOPT \cite{Acar_SDM_2010}, the state-of-the-art batch algorithm. The relative change in function value tolerance is set to $10^{-9}$ and the maximum iterations is 300 for CP-WOPT. Our proposed algorithm gives superior convergence performances than those of TeCPSGD.

\begin{figure*}[htbp]
\centering
\begin{minipage}[b]{.49\linewidth}
\centering
\centerline{\includegraphics[width=1.0\linewidth]{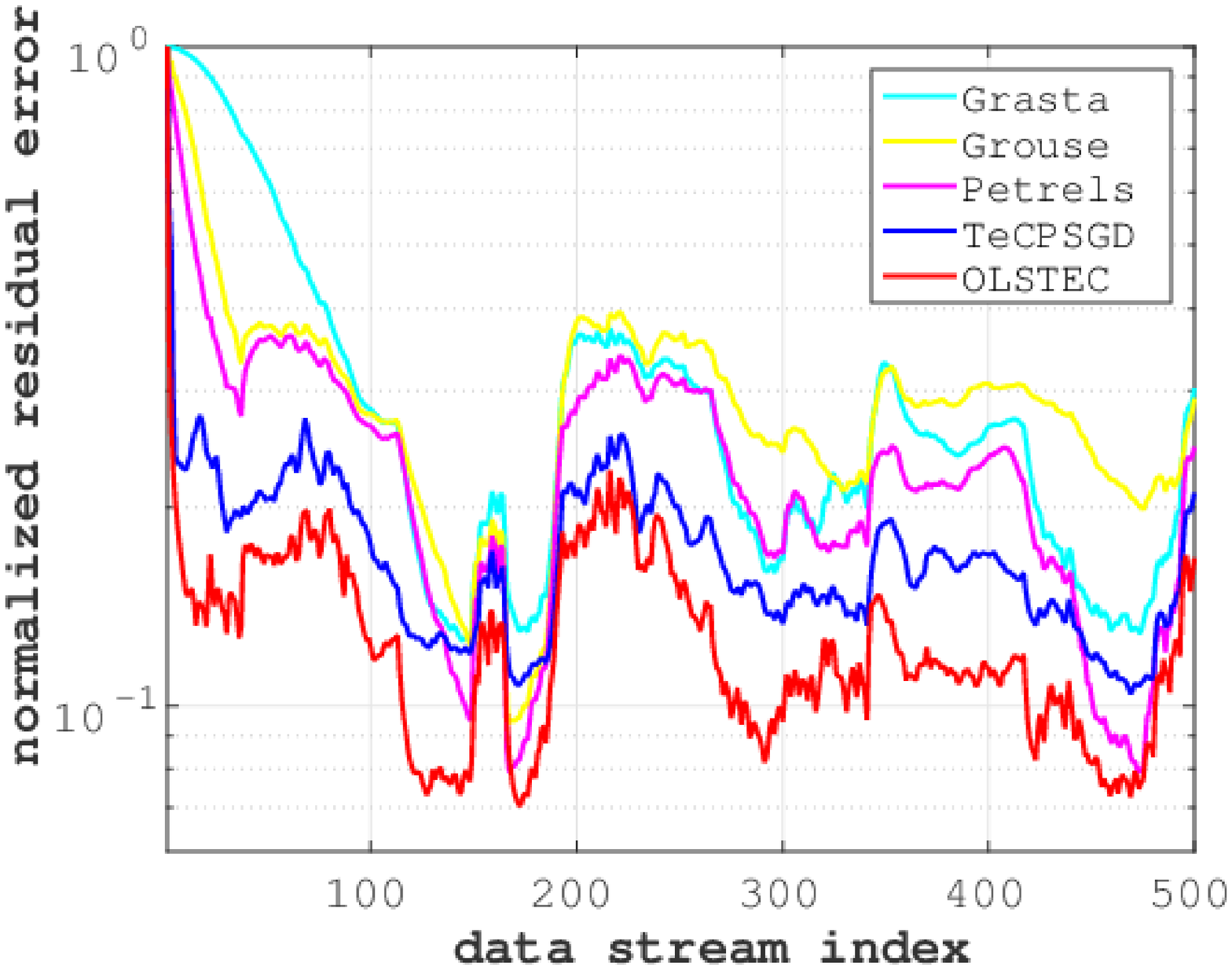}}
\vspace*{-1.5cm}
\centerline{\scriptsize (a) Stationary background ($\rho = 0.1$)}\medskip
\end{minipage}
\hspace{-0.2cm}
\begin{minipage}[b]{.49\linewidth}
\centering
\centerline{\includegraphics[width=1.0\linewidth]{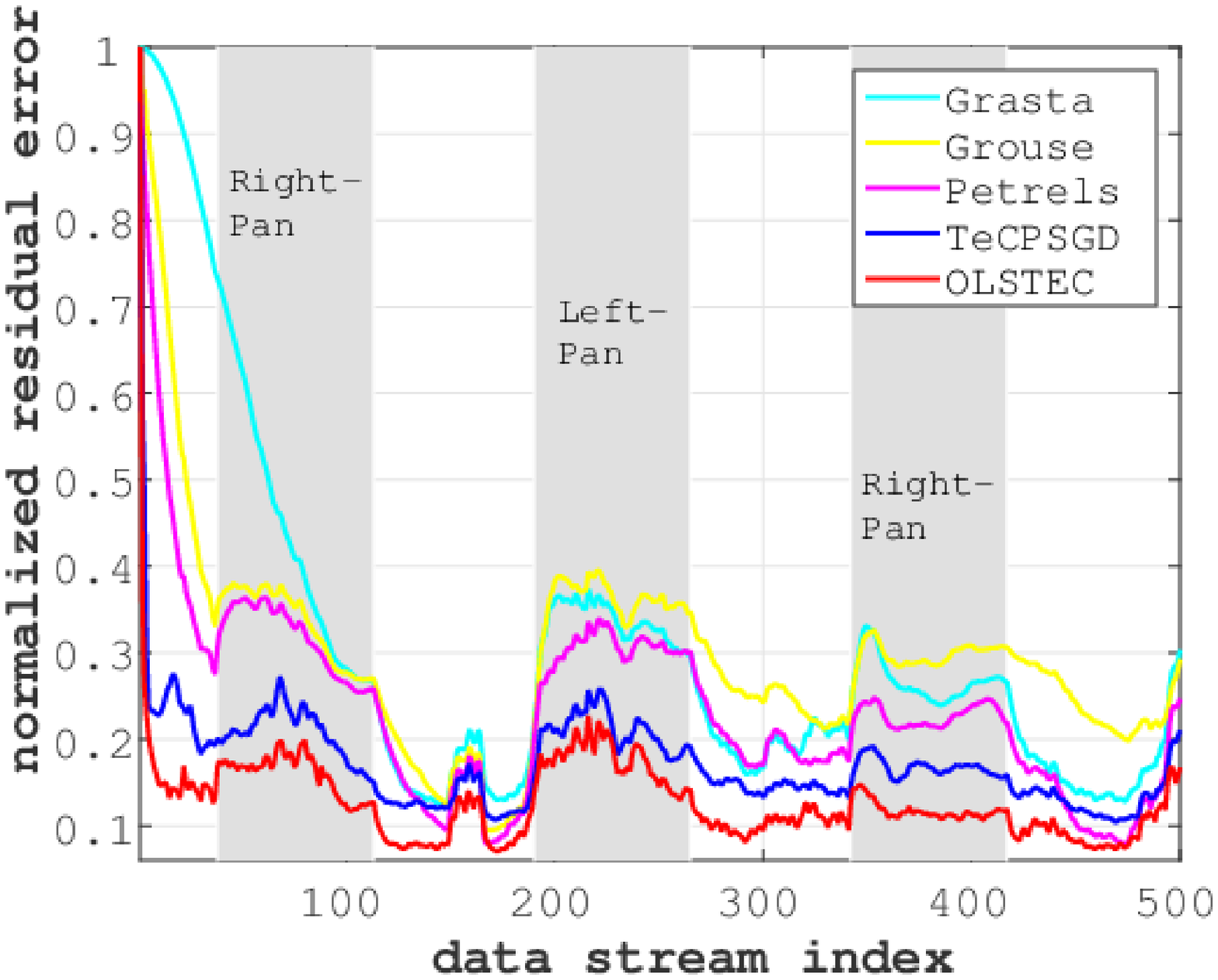}}
\vspace*{-1.5cm}
\centerline{\scriptsize (b) Dynamic background ($\rho = 0.1$)}\medskip
\end{minipage}\\
\vspace*{1.0cm}

\begin{minipage}[b]{0.6\linewidth}
\centerline{\vspace*{2cm}\includegraphics[width=0.9\linewidth]{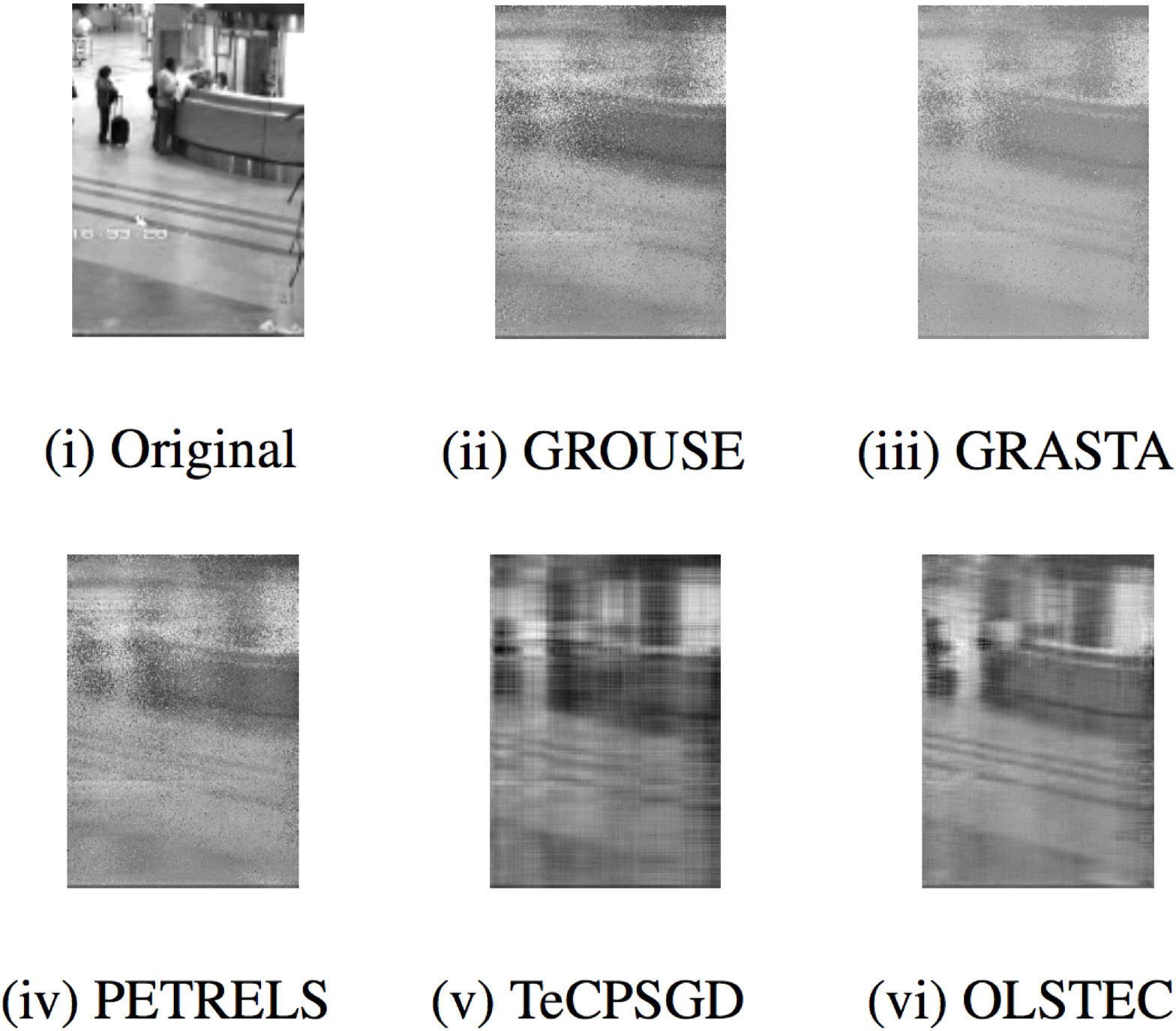}}
\vspace*{-1.8cm}
 \centerline{\scriptsize (c) Reconstructed subspace images.}\medskip
\end{minipage}
\caption{The normalized estimation error in real-world dataset.}
\label{fig:realworld_data}
\end{figure*}

We also evaluate a scenario where a subspace of rank changes abruptly periodically. Four rank-5 tensors of $100\times100\times250$ are concatenated in series at the 3-rd order direction. Figure \ref{fig:synthetic_data}(c) shows the normalized residual error at each iteration. This shows that the subspace tracking behavior of the OLSTEC algorithm gives a superior performance than that of TeCPSGD which cannot recover correct subspaces after abrupt changes.

Next, we evaluate the tracking performances using surveillance video as a real-world dataset.  Although each video frame does not have low-rank structure and a tensor-based approach basically has a disadvantage for the approximation of its underlying subspace, this experiments demonstrates the superior tacking performance of OLSTEC. We compare OLSTEC with TeCPSGD as well as the matrix-based algorithms including GROUSE \cite{Balzano_arXiv_2010_s}, GRASTA \cite{He_CVPR_2012_s}, and PETRELS \cite{Chi_IEEETransSP_2013}. We use Matlab codes provided by the respective authors except for TeCPSGD with our configured parameters. ``Airport Hall" dataset of size $288\times 352$ with $500$ frames is used. Moreover, for fair comparison between tensor and matrix-based algorithms, the rank is set to $20$ and $10$ for the former, i.e., OLSTEC and TeCPSGD, and for the latter, respectively. Still, the tensor-based algorithms has much less free parameters than those of the matrix-based algorithms. 
This experiment also considers two scenarios.
The first separates foreground objects with static background and moving objects in the foreground. Figure \ref{fig:realworld_data} (a) shows the superior performance of OLSTEC against other algorithms. 
Furthermore, we examine the performances against a dynamic moving background as the second scenario. The input video is created virtually by moving cropped partial image from its original entire frame image of video. The cropping window with $288 \times 200$ moves from the leftmost partial image to the rightmost, then returns to the leftmost image after stopping a certain period of time. 
The generated video includes right-panning video from $38$-th to $113$-th frame and from $342$-th to $417$-th frame, and left-panning video from $190$-th to $265$-th frame.
Figure \ref{fig:realworld_data}(b) shows how OLSTEC can quickly adapt to the changed background.
Figure \ref{fig:realworld_data}(c) shows the reconstructed (i.e., completed) image at $110$-th frame  of OLSTEC gives better quality than those of others.

\section{Conclusion and future work}

We have proposed a new online tensor subspace tracking algorithm, dubbed OLSTEC, for the partially observed high-dimensional data stream corrupted by noise. Especially, we addressed a second-order stochastic gradient descent based on the recursive least squares to achieve faster convergence of subspace tracking. Numerical comparisons suggest that our proposed algorithm has superior performances on synthetic as well as real-world datasets.  As a future research direction, we will investigate the ways of the Tucker decomposition. 

\vfill\pagebreak

\section*{Acknowledgments}
H. Kasai thanks Prof. Wolfgang Kellerer and Prof. Martin Kleinsteuber for useful discussions on the paper. H. Kasai is (partly) supported by the Ministry of Internal Affairs and Communications, Japan, as the SCOPE Project.

\bibliographystyle{IEEEtran}
\bibliography{kasai}

\begin{thebibliography}{10}
\providecommand{\url}[1]{#1}
\csname url@samestyle\endcsname
\providecommand{\newblock}{\relax}
\providecommand{\bibinfo}[2]{#2}
\providecommand{\BIBentrySTDinterwordspacing}{\spaceskip=0pt\relax}
\providecommand{\BIBentryALTinterwordstretchfactor}{4}
\providecommand{\BIBentryALTinterwordspacing}{\spaceskip=\fontdimen2\font plus
\BIBentryALTinterwordstretchfactor\fontdimen3\font minus
  \fontdimen4\font\relax}
\providecommand{\BIBforeignlanguage}[2]{{%
\expandafter\ifx\csname l@#1\endcsname\relax
\typeout{** WARNING: IEEEtran.bst: No hyphenation pattern has been}%
\typeout{** loaded for the language `#1'. Using the pattern for}%
\typeout{** the default language instead.}%
\else
\language=\csname l@#1\endcsname
\fi
#2}}
\providecommand{\BIBdecl}{\relax}
\BIBdecl

\bibitem{Liu_IEEETransPAMI_2013_s}
J.~Liu, P.~Musialski, P.~Wonka, and J.~Ye, ``Tensor completion for estimating
  missing values in visual data,'' \emph{IEEE Trans. Pattern Anal. Mach.
  Intell.}, vol.~35, no.~1, pp. 208--220, 2013.

\bibitem{Tomioka_Latent_2011_s}
R.~Tomioka, K.~Hayashi, and H.~Kashima, ``Estimation of low-rank tensors via
  convex optimization,'' \emph{arXiv:1010.0789}, 2011.

\bibitem{Signoretto_MachineLearning_2014_s}
M.~Signoretto, Q.~T. Dinh, L.~D. Lathauwer, and J.~A. Suykens, ``Learning with
  tensors: a framework based on convex optimization and spectral
  regularization,'' \emph{Mach. Learn.}, vol.~94, no.~3, pp. 303--351, 2014.

\bibitem{Candes_FoundCompuMath_2009_s}
E.~J. Cand{\`e}s and B.~Recht, ``Exact matrix completion via convex
  optimization,'' \emph{Found. Comput. Math.}, vol.~9, no.~6, pp. 717--772,
  2009.

\bibitem{Filipovi_MultiSysSigPro_2013_s}
M.~Filipovi{\'c} and A.~Juki{\'c}, ``{Tucker} factorization with missing data
  with application to low- n -rank tensor completion,'' \emph{Multidim. Syst.
  Sign. P.}, 2013.

\bibitem{Kressner_BIT_2014_s}
D.~Kressner, M.~Steinlechner, and B.~Vandereycken, ``Low-rank tensor completion
  by {Riemannian} optimization,'' \emph{BIT Numer. Math.}, vol.~54, no.~2, pp.
  447--468, 2014.

\bibitem{Boumal_NIPS_2011}
N.~Boumal and P.-A. Absil, ``{RTRMC} : A {Riemannian} trust-region method for
  low-rank matrix completion,,'' in \emph{Proceedings of the Annual Conference
  on Neural Information Processing Systems (NIPS)}, 2011.

\bibitem{Mishra_SIAMOpt_2013}
B.~Mishra, G.~Meyer, F.~Bach, and R.~Sepulchre, ``Low-rank optimization with
  trace norm penalty,'' \emph{SIAM Journal on Optimization}, vol.~23, no.~4,
  pp. 2124--2149, 2013.

\bibitem{Ngo_NIPS_2012_s}
T.~Ngo and Y.~Saad, ``Scaled gradients on {Grassmann} manifolds for matrix
  completion,'' in \emph{NIPS}, 2012, pp. 1421--1429.

\bibitem{Yang_IEEESP_1995_s}
B.~Yang, ``Projection approximation subspace tracking,'' \emph{IEEE Trans. on
  Signal Processing}, vol.~43, no.~1, pp. 95--107, 1995.

\bibitem{Balzano_arXiv_2010_s}
L.~Balzano, R.~Nowak, and B.~Recht, ``Online identification and tracking of
  subspaces from highly incomplete information,'' \emph{arXiv:1006.4046}, 2010.

\bibitem{Edelman98a}
A.~Edelman, T.~Arias, and S.~Smith, ``The geometry of algorithms with
  orthogonality constraints,'' \emph{SIAM J. Matrix Anal. Appl.}, vol.~20,
  no.~2, pp. 303--353, 1998.

\bibitem{Absil_OptAlgMatManifold_2008}
P.-A. Absil, R.~Mahony, and R.~Sepulchre, \emph{Optimization Algorithms on
  Matrix Manifolds}.\hskip 1em plus 0.5em minus 0.4em\relax Princeton
  University Press, 2008.

\bibitem{He_CVPR_2012_s}
J.~He, L.~Balzano, and A.~Szlam, ``Incremental gradient on the grassmannian for
  online foreground and background separation in subsampled video,'' \emph{IEEE
  Conference on Computer Vision and Pattern Recognition (CVPR)}, 2012.

\bibitem{Chi_IEEETransSP_2013}
Y.~Chi, Y.~C. Eldar, and R.~Calderbank, ``Petrels: Parallel subspace estimation
  and tracking using recursive least squares from partial observations,''
  \emph{IEEE Trans. on Signal Processing}, vol.~61, no.~23, pp. 5947--5959,
  2013.

\bibitem{Nion_IEEETransSP_2009}
D.~Nion and N.~Sidiropoulos, ``Adaptive algorithms to track the parafac
  decomposition of a third-order tensor,'' \emph{IEEE Transactions on Signal
  Processing}, vol.~57, no.~6, pp. 2299--2310, 2009.

\bibitem{Yu_ICML_2015_s}
R.~Yu, D.~Cheng, and Y.~Liu, ``Accelerated online low-rank tensor learning for
  multivariate spatio-temporal streams,'' \emph{International Conference on
  Machine Learning (ICML)}, 2015.

\bibitem{Mardani_IEEETransSP_2015}
M.~Mardani, G.~Mateos, and G.~Giannakis, ``Subspace learning and imputation for
  streaming big data matrices and tensors,'' \emph{IEEE Transactions on Signal
  Processing}, vol.~63, no.~10, pp. 266--2677, 2015.

\bibitem{Kolda_SIAMReview_2009_s}
T.~G. Kolda and B.~W. Bader, ``Tensor decompositions and applications,''
  \emph{SIAM Review}, vol.~51, no.~3, pp. 455--500, 2009.

\bibitem{Acar_SDM_2010}
E.~Acar, D.~M. Dunlavy, T.~G. Kolda, and M.~M{\o}rup, ``Scalable tensor
  factorizations with missing data,'' in \emph{Proceedings of the 2010 SIAM
  International Conference on Data Mining (SDM10)}, 2010, pp. 701--712.

\end{thebibliography}

\end{document}